\documentclass[12pt]{article}
\usepackage{hyperref}
\usepackage{amsfonts}
\usepackage{enumerate}
\begin{document}

\begin{center}
{\large\bf A note on the causality of singular linear discrete time systems}

\vskip.20in
Christos Tsegkis$^{1}$\\[2mm]
{\footnotesize
$^{1}$ School of Engineering, The University of Edinburgh, UK}
\end{center}

{\footnotesize
\noindent
\textbf{Abstract:}
In this article we study the causality of non-homogeneous linear singular discrete time systems whose coefficients are square constant matrices. By assuming that the input vector changes only at equally space sampling instants we provide properties for causality between state and inputs and causality between output and inputs.
\\
\\[3pt]
\\[3pt]

\vskip.2in

\section{Introduction}

In this article we shall be concerned with the non-homogeneous singular discrete time system of the form
\begin{equation}
\begin{array}{c}FY_{k+1}=GY_k+BV_k\\
X_k=CY_k\end{array}
\end{equation}
with known initial conditions
\begin{equation}
Y_{k_0}
\end{equation}        
where $F,G \in \mathcal{M}(n \times n;\mathcal{F})$, $Y_k
\in \mathcal{M}(n \times 1;\mathcal{F})$ (i.e., the algebra of
square matrices with elements in the field $\mathcal{F}$), $X_k
\in \mathcal{M}(m \times 1;\mathcal{F})$, $B \in \mathcal{M}(n \times l;\mathcal{F})$ and $C \in \mathcal{M}(m \times n;\mathcal{F})$. For the
sake of simplicity, we set ${\mathcal{M}}_n  = {\mathcal{M}}({n
\times n;\mathcal{F}})$ and ${\mathcal{M}}_{nm}  =
{\mathcal{M}}({n \times m;\mathcal{F}})$. We assume that the system (1) is singular, i.e. the matrix $F$ is singular and that the input vector $V_k$ changes only at equally space sampling instants. Many authors have studied discrete time systems, see and their applications, see [1-9, 12, 13, 16-18, 20-32]. In this article we study the causality of these systems. The results of this paper can be applied also in systems of fractional nabla difference equations, see [10, 11]. In addition they are very useful for applications in many mathematical models using systems of difference equations existing in the literature, see [14, 15, 29-32]. 
\\\\
\textbf{Definition 1.1.} Given $F, G \in \mathcal{M}_{nm}$  and an indeterminate $s\in \textsl{F}$, the matrix
pencil $sF-G$ is called regular when $m=n$ and  $\det(sF-G)\neq 0$.
In any other case, the pencil will be called singular.
\\\\
In this article, we consider the case that the pencil is \emph{regular}. The class of the pencil $sF-G$ is characterized by a uniquely
defined element, known as a complex Weierstrass canonical form,
$sF_w -Q_w$, see [19, 24], specified by the complete set of
invariants of $sF-G$. This is the set of \emph{elementary divisors} (e.d.). In the case
of a regular matrix pencil, we have e.d. of the following type:
\begin{itemize}
    \item e.d. of the type  $(s-a_j)^{p_j}$, \emph{are called finite elementary
    divisors} (f.e.d.), where $a_j$ is a finite eigenavalue of algebraic multiplicity $p_j$
    
    \item e.d. of the type  $\hat{s}^q$ are called \emph{infinite elementary divisors}
    (i.e.d.), where q the algebraic multiplicity of the infinite eigenvalues
\end{itemize}
We assume that $\sum_{i =1}^\nu p_j  = p$ and $p+q=n$.\\ Let $B_1 ,B_2 ,\dots, B_n $ be elements of $\mathcal{M}_n$. The direct sum
of them denoted by $B_1  \oplus B_2  \oplus \dots \oplus B_n$ is
the blockdiag$\left[\begin{array}{cccc} B_1& B_2& \dots& B_n\end{array}\right]$ . From the regularity of $sF-G$, there exist nonsingular matrices $P, Q\in \mathcal{M}_n$ such that 
\[
PFQ = F_w  = I_p  \oplus H_q
\]
and
\[
PGQ = G_w  = J_p  \oplus I_q
\]
Where $sF_w -Q_w$ is the complex Weierstrass form of the regular
pencil sF-G and is defined by $sF_w  - Q_w :=sI_p  - J_p  \oplus sH_q  - I_q $, where the first normal Jordan
type element is uniquely defined by the set of the finite eigenvalues.
\[
  ({s - a_1 })^{p_1 } , \dots ,({s - a_\nu  }
 )^{p_\nu }
\]
of $sF-G$. The second block has the form
\[
    sI_p  - J_p  := sI_{p_1 }  - J_{p_1 } (
    {a_1 }) \oplus  \dots  \oplus sI_{p_\nu  }  - J_{p_\nu  }
    ({a_\nu  }) 
\]
And also the q blocks of the third uniquely defined block
$sH_q -I_q$ correspond to the infinite eigenvalues
\[
  \hat s^{q_1} , \dots ,\hat s^{q_\sigma}, \quad \sum_{j =
  1}^\sigma  {q_j  = q}
\]
of $sF-G$  and has the form
\[
    sH_q  - I_q  := sH_{q_1 }  - I_{q_1 }  \oplus
    \dots  \oplus sH_{q_\sigma  }  - I_{q_\sigma}
\]
Thus, $H_q$  is a nilpotent element of $\mathcal{M}_n$  with index
$\tilde q = \max \{ {q_j :j = 1,2, \ldots ,\sigma } \}$, where
\[
    H^{\tilde q}_q=0_{q, q},
\]
and $J_{p_j } ({a_j }),H_{q_j }$ are defined as

   \[
   J_{p_j } ({a_j }) =  \left[\begin{array}{ccccc}
   a_j  & 1 & \dots&0  & 0  \\
   0 & a_j  &   \dots&0  & 0  \\
    \vdots  &  \vdots  &  \ddots  &  \vdots  &  \vdots   \\
   0 & 0 &  \ldots& a_j& 1\\
   0 & 0 & \ldots& 0& a_j
   \end{array}\right] \in {\mathcal{M}}_{p_j }, H_{q_j }  = \left[
\begin{array}{ccccc} 0&1&\ldots&0&0\\0&0&\ldots&0&0\\\vdots&\vdots&\ddots&\vdots&\vdots\\0&0&\ldots&0&1\\0&0&\ldots&0&0
\end{array}
\right] \in {\mathcal{M}}_{q_j }.
\]
For algorithms about the computations of the jordan matrices see [19, 24, 28].

\section{The solution of a singular linear discrete time system}

 In this subsection we obtain formulas for the solutions of LMDEs with regular matrix pencil and we
give necessary and sufficient conditions for existence and uniqueness of solutions.
\\\\
\textbf{Theorem 2.1.}
Consider the system (1), (2) and  let the $p$ linear independent (generalized) eigenvectors of the finite eigenvalues of the pencil sF-G be the columns of a matrix $Q_p$. Then the solution is unique if and only if
\begin{equation}
Y_{k_0}\in colspan Q_p+QD_{k_0}
\end{equation}
Moreover the analytic solution is given by
\begin{equation}
    Y_k=Q_pJ_p^{k-k_0}Z^p_{k_0}+QD_k
\end{equation}
where $D_k=\left[
\begin{array}{c} \sum^{k-1}_{i=0}J_p^{k-i-1}B_pV_i\\-\sum^{q_{*}-1}_{i=0}H_q^iB_qV_{k+i}
\end{array}\right]$ and PB=$\left[\begin{array}{c} B_p \\ B_q\end{array}\right]$, with $B_p \in \mathcal{M}_{pn}$, $B_q \in \mathcal{M}_{qn}$.
\\\\
\textbf{Proof.} Consider the transformation
\[
    Y_k=QZ_k
\]
Substituting the previous expression into (1) we obtain
\[
    FQZ_{k+1}=GQZ_k+BV_k.
\]
Whereby, multiplying by  P, we arrive at
\[
    F_wZ_{k+1}=G_w Z_k+PBV_k.
\]
Moreover, we can write $Z_k$ as
$Z_k =\left[\begin{array}{c}
 Z^p_k \\
 Z_K^q
 \end{array}\right] $.
Taking into account the above expressions, we arrive easily at two subsystems of (1). The subsystem
\begin{equation}
    Z_{k+1}^p = J_p Z_k^p+B_pV_k 
 \end{equation}
and the subsystem
\begin{equation}
    H_q Z_{k+1}^q = Z_k^q+B_q V_k
\end{equation}
The subsystem (5) has the unique solution 
\begin{equation}
    Z_k^p=J_p^{k-k_0}Z^p_{k_0}+\sum^{k-1}_{i=0}J_p^{k-i-1}B_pV_i , k\geq k_0, 
\end{equation}
see [1, 4, 11, 12]). By applying the Zeta transform we get the solution of subsystem (6)  
\begin{equation}
Z_k^q=-\sum^{q_*-1}_{i=0}H_q^iB_qV_{k+i}
\end{equation}
 Let $Q=\left[\begin{array}{cc} Q_p& Q_q\end{array} \right]$, where $Q_p \in \mathcal{M}_{np}$, $Q_q\in \mathcal{M}_{nq}$ the matrices with columns the p, q generalized eigenvectors of the finite and infinite eigenvalues respectively. Then we obtain
\[
     Y_k = QZ_k =
     [Q_p  Q_q ]
     \left[\begin{array}{c}
     J_p^{k-k_0}Z^p_{k_0}+\sum^{k-1}_{i=0}J_p^{k-i-1}B_pV_i  \\
     -\sum^{q_*-1}_{i=0}H_q^iB_qV_{k+i}
    \end{array}\right]  
    \]
    \[
    Y_k=Q_pJ_p^{k-k_0}Z^p_{k_0}+Q_p\sum^{k-1}_{i=0}J_p^{k-i-1}B_pV_i-Q_q-\sum^{q_*-1}_{i=0}H_q^iB_qV_{k+i} .
    \]
    \[
    Y_k=Q_pJ_p^{k-k_0}Z^p_{k_0}+QD_k
    \]
The solution that exists if and only if 
\[
Y_{k_0}=Q_pZ_{k_0}^p+QD_{k_0}
\]
or
\[
Y_{k_0}\in colspan Q_p+QD_{k_0}
\]

\section{Causality}

Generally for systems of type (1) we define the notion of causality, which is properly defined bellow
\\\\
\textbf{Definition 3.1.} The non-homogeneous singular continuous system (1) is called casual, if its state $Y_k$, for any $k > k_0$ is determined completely by initial state $Y_{k_0}$ and former inputs $V_{k_0}$, $V_{k_0+1}$, ..., $V_k$. Otherwise it is called noncausal.
\\\\
Discrete time normal systems are characterized by the property of causality. Next we will study the causality in a singular system of the form (1).
\subsection*{Causality between state and inputs}
\textbf{Proposition 3.1.} In system (1) causality between state and inputs exists if and only if $H_qB_q=0_{q, l}$
\\\\
\textbf{Proof.} From (8) it is clear that the state $Z_k$ and obviously $Y_k$ for any $k\geq k_0$ is to be determined by former inputs if and only if $H_q^iB_q=0_{q,l}$ for every $i=1, 2, ..., q^*-1$, which is equivalent to the relation $H_qB_q=0_{q, l}$.

\subsection*{Causality between output and inputs}
\textbf{Proposition 3.2.} In system (1) causality between output and input exists if and only if 
\begin{equation}
CQ_qH_q^iB_q=0_{m,l}
\end{equation}
for every $i= 1, 2, ..., q^*-1$.
\\\\
\textbf{Proof.} The solution of the state equation of the system (1) is given by Theorem 2.1. 
\[
Y_k =Q_pJ_p^{k-k_0}Z^p_{k_0}+Q_p\sum^{k-1}_{i=0}J_p^{k-i-1}B_pV_i-Q_q\sum^{q_*-1}_{i=0}H_q^iB_qV_{k+i}
\]
Setting the expression of $Y_k$ in the state output relation $X_k=CY_k$ we take
\begin{equation}
X_k =CQ_pJ_p^{k-k_0}Z^p_{k_0}+CQ_p\sum^{k-1}_{i=0}J_p^{k-i-1}B_pV_i-CQ_q\sum^{q_*-1}_{i=0}H_q^iB_qV_{k+i}
\end{equation}
From the above expression it is clear that non causality is due to the existence of the term $\sum^{q_*-1}_{i=0}CQ_qH_q^iB_qV_{k+i}$. So the causal relationship between $X_k$ and $V_k$ exists if and only if $CQ_qH_q^iB_q=0_{m,l}$ for every $i= 1, 2, ..., q^*-1$.
\\\\
The relation (9) can be written equivalently as 
\begin{equation}
C\left[\begin{array}{ccc}Q_qH_qB_q&...&Q_qH_q^{q^*-1}B_q\end{array}\right]=0_{m,q^*nl}
\end{equation}
So the following Proposition is obvious.
\\\\
\textbf{Proposition 3.3.} The system (1) is causal if and only if every column of the matrix $\left[\begin{array}{ccc}Q_qH_qB_q&...&Q_qH_q^{q^*-1}B_q\end{array}\right]$ lies in the right nullspace of the matrix $C$.
\\\\
\textbf{Remark 3.1.} If the system pencil $sF-G$ has no infinite eigenvalues then the matrix $Q_q=0_{n,q}$. So the relation (11) is satisfied and we have causality between inputs and outputs of the system.

\section*{Conclusions}
Having shown that the solution of the discrete time system of the form (1) exists if the initial conditions (2) belong to the set (3) and is given by the formula (4), we prove that in system (1) causality between state and inputs and causality between output and inputs exists under necessary and sufficient conditions.

\subsection*{Acknowledgments} I would like to express my
sincere gratitude to Professor G.I. Kalogeropoulos and Dr. I.K. Dassios for their helpful and
fruitful discussions that improved this article.

\end{document}